\documentclass{amsart}

\usepackage[utf8]{inputenc}
\usepackage[english]{babel}
\usepackage{hyperref}

\usepackage{amssymb, amsfonts, amsmath}
\usepackage{graphicx}
\usepackage{float, caption}
\usepackage{amsthm}

\title{Inverse-Magnetic Billiards on a Square}

\author{Andres Perico}
\email{aperico@ucsc.edu}
\address{Department of Mathematics, University of California, Santa Cruz, CA 95060}

\begin{document}
	\maketitle
\section{Introduction}
	We work in a unit square $\Omega\subset \mathbb{R}^2$ where we have a constant magnetic field outside of $\Omega$ with magnitude $B$, inside there is no magnetic field. A particle at an initial position in the boundary of $\Omega$ starts moving towards the interior of the square. Without loss of generality we can assume that the square is in canonical position with vertices $(0,0),(1,0),(1,1),(0,1)$ and that the particle starts at some point on the bottom side (on the $x-axis$).
	
	\begin{minipage}{\linewidth}
	\makebox[\linewidth]{
		\includegraphics[width=4in]{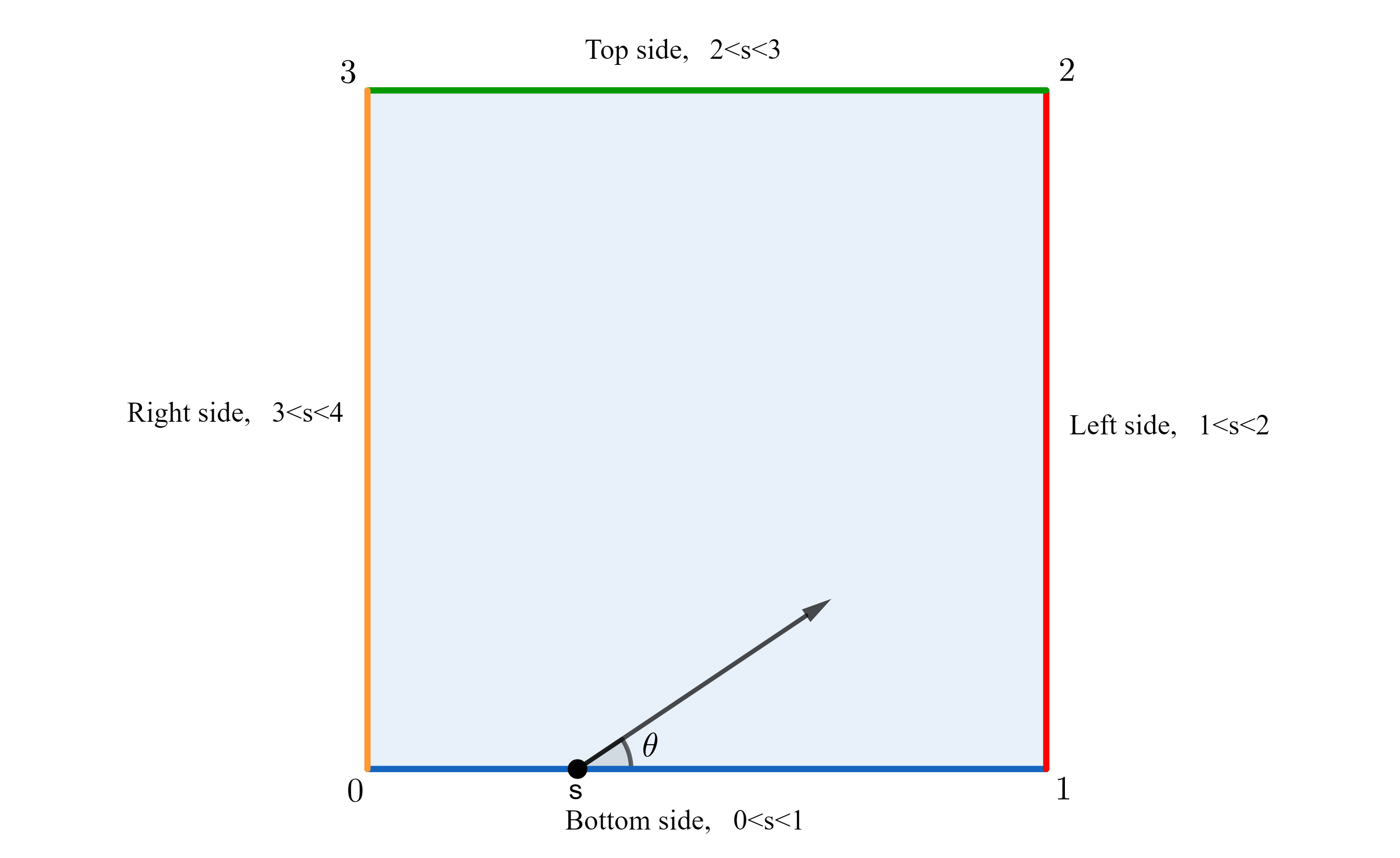}}
		\captionof{figure}{Initial conditions}
	\end{minipage}
	
	\vspace{10pt}
	\subsection{Description of dynamics}
	Say that the initial conditions are $(s,\theta)$ where $s\in (0,1)$, $\theta\in (0,\pi)$ (we are going to omit the corners for now). Following the billiards notation, these coordinates for initial conditions are called Birkhoff coordinates. The particle moves inside $\Omega$ in straight line at an angle $\theta$ with respect to the boundary until it hits the boundary again at $(s_1,\theta_1)$ on the \textit{exiting} side. It has three options for exiting side, the other three sides of the square.
	We are considering an electron with charge $-1$, this is why the circular motion outside of $\Omega$ corresponds to a counterclockwise movement.
	After hitting the boundary it moves on a circle of radius $r=1/B$ tangent at the exit point to the line of its previous trajectory.\\
	The particle will hit $\Omega$ again at $(s_2,\theta_2)$ on the \textit{entering} side, depending on $B$ the radius of rotation will be smaller or bigger giving different return points on different sides. The map $(s_1,\theta_1)\to (s_2,\theta_2)$ will be referred to as the \textit{magnetic bounce}.
	As soon as it hits $\Omega$ again the particle moves in straight line, this line is tangent to the circle of its previous trajectory at the point of entry (see figure \ref{iteration}).
	
	\begin{minipage}{\linewidth}
	\makebox[\linewidth]{
		\includegraphics[width=4in]{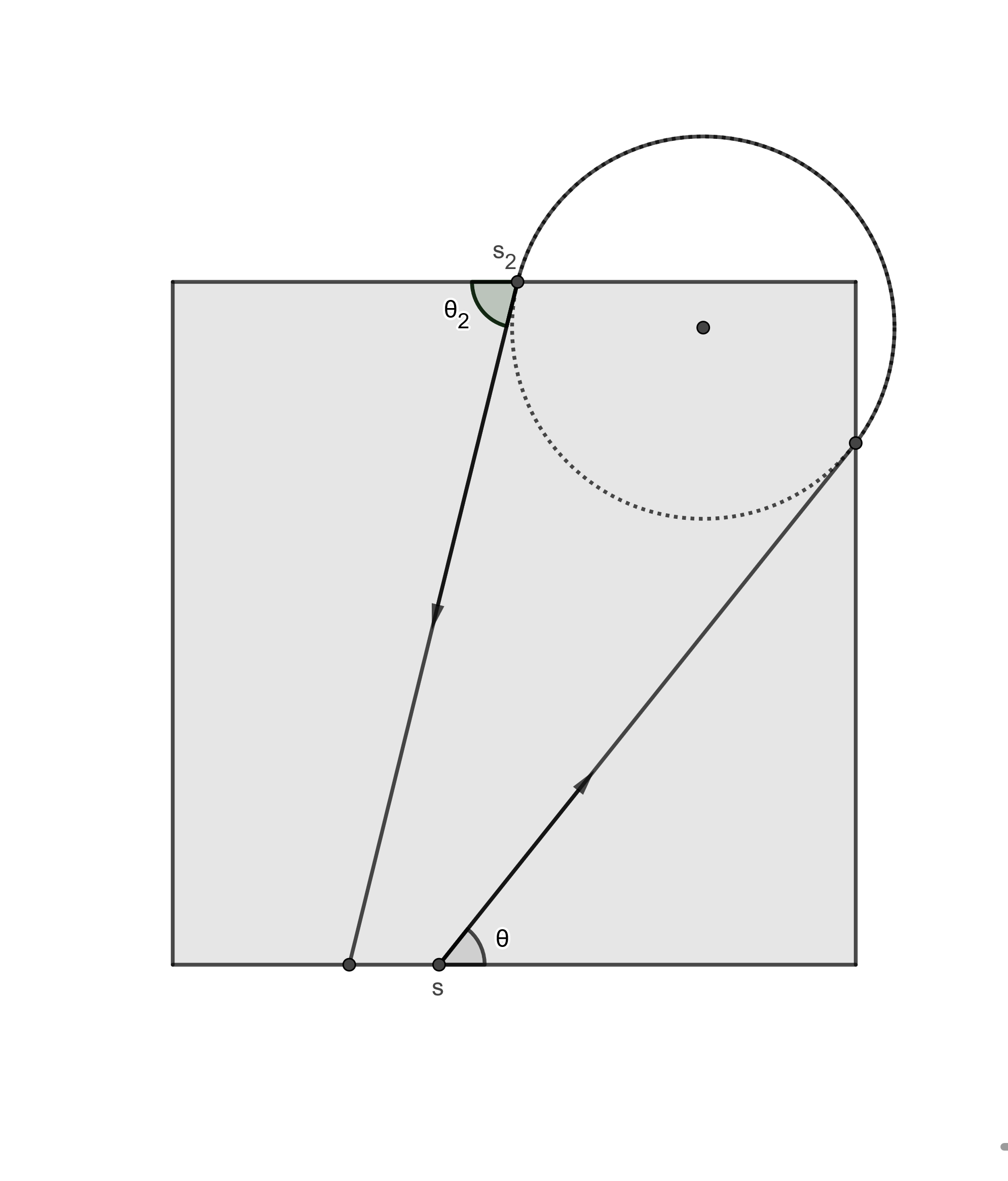}}
		\vspace{-50pt}
		\captionof{figure}{One bounce in the magnetic billiard}\label{iteration}
	\end{minipage}
	
	\vspace{10pt}
	
	This dynamics can be described by a map
	\[ F: \Sigma^2 \rightarrow \Sigma^2 \]
	\[ F(s,\theta)=(s_2,\theta_2) \]
	where $\Sigma^2=\{ (s,\theta) | 0\leq s \leq 4, \quad 0<\theta <\pi \} $ is the space of directed unit vectors toward the interior of $\Omega$ with initial points on the boundary.
	
	We want to study the dynamics of this billiard type map. 
	The most challenging aspect of this investigation concerns trajectories which go around a corner.  See figure \ref{iteration}. The map is apparently chaotic and all the chaos seems to come from turning around corners.\\
	
	Here are some of the questions we want to answer:
	
	\begin{itemize}
		\item Are there any periodic orbits?
		\item For any $B$, is there a periodic orbit?
		\item Can we classify the periodic orbits? 
		\item Can we classify non periodic orbits?
		\item Is any orbit dense in $\Omega$?
		\item Is any orbit dense in $\Sigma^2$?
		\item As a map $F:\Sigma^2\rightarrow \Sigma^2$ symplectic integrable.  
		\item How do orbits look like?
	\end{itemize}

\subsection{Extreme cases: $B\to 0,\infty$.\\}	
	Classic billiards is the limit case $\displaystyle r= \lim_{B \to \infty} r_B$, so the particle is trapped in the square and bounces in the classic way.\\
	The case when $\displaystyle r =\lim_{B \to 0} r_B$ refers when the particle continues an straight trajectory all the time, so the only part of the trajectory inside the square is the initial segment. This orbit has only one exit and one entry point: the entry point is exactly the same starting point, the particle get there again after infinite time.
	
\section{No turning corners}
	
	We are assuming our initial condition is on the bottom side of the square ($0<s<1$), we look for the returning map $F$ depending on the exiting side: right, top or left side. Assume that the entry point is on the same side as the exit point, this means that the trajectory outside $\Omega$ doesn't go around a corner.\\
	
	Figure \ref{slope1nocorners} shows an example of this case, $\theta = \frac{\pi}{4}$ or initial slope equal to $1$ and radius small enough (depending on our initial conditions) so the magnetic bounce happens on the same side of the square.
	
	\begin{minipage}{\linewidth}
	\makebox[\linewidth]{	
		\includegraphics[width=5in]{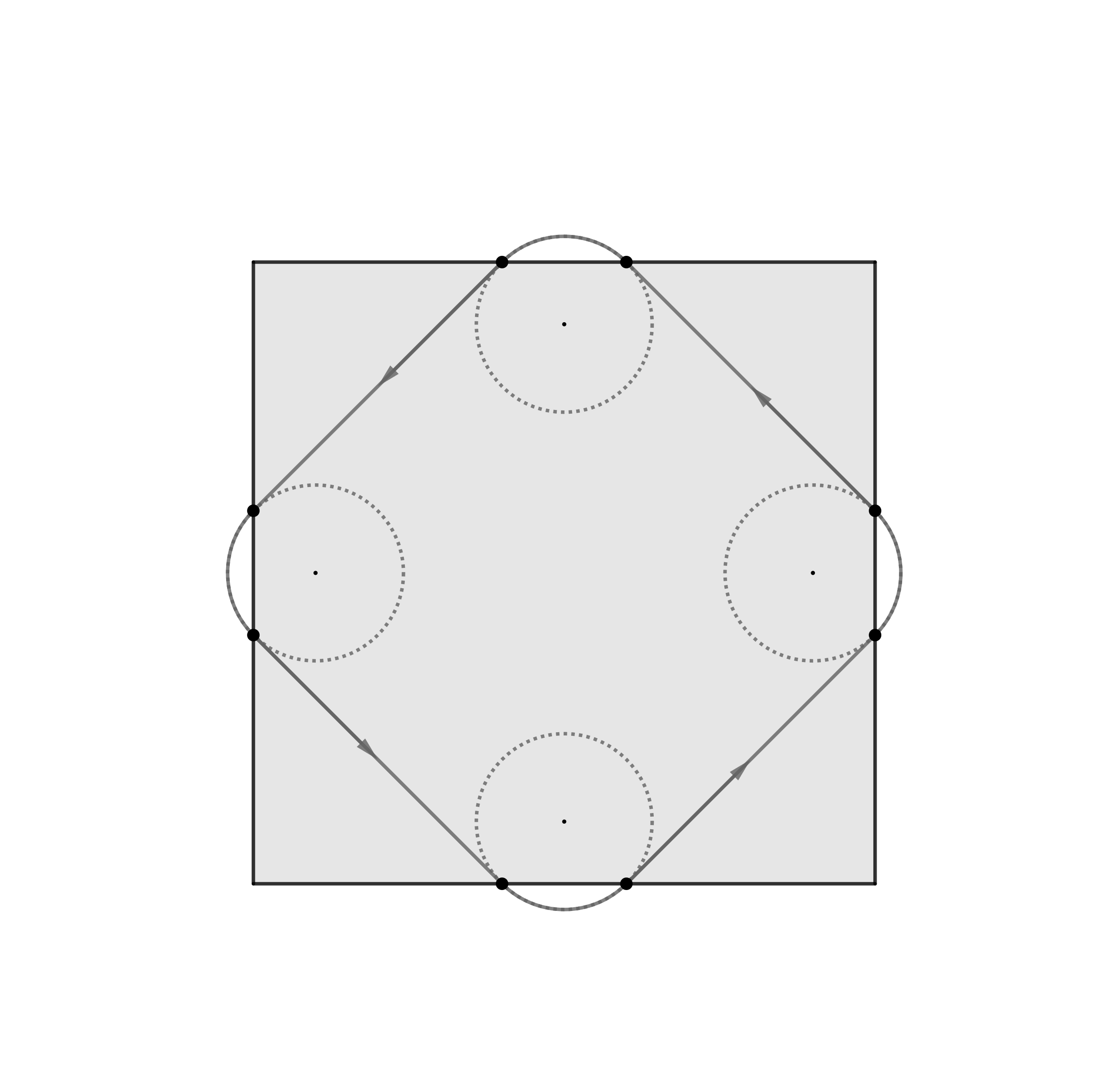}}
		\vspace{-50pt}
		\captionof{figure}{Example of no turning corners}\label{slope1nocorners}
	\end{minipage}
	
	\newpage

	\begin{minipage}{\textwidth}
		\begin{center}
			\bf{Bouncing on the right side of the square}
		\end{center}
	\makebox[\linewidth]{	
		\includegraphics[width=3.6in]{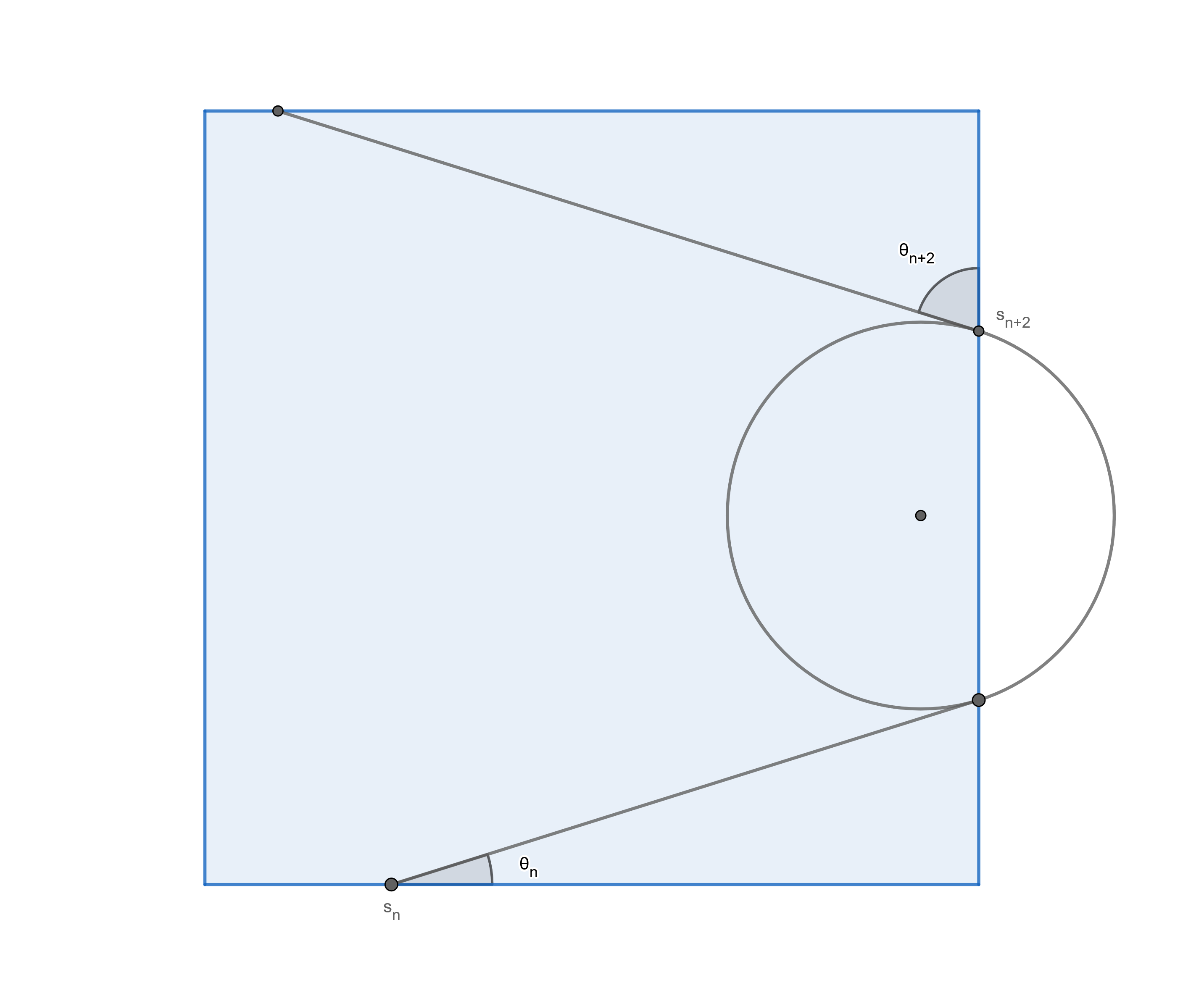}}
		\vspace{-20pt}
		\captionof{figure}{Bouncing on the right side}\label{rightside}
	\begin{equation*}
	s_{n+1} = (1-s_n)\tan\theta, \qquad
	s_{n+2} = s_{n+1} + \frac{2}{B}\cos\theta,\qquad
	\theta_{n+2} = \frac{\pi}{2}-\theta_n.
	\end{equation*}
	\end{minipage}

	\vspace{5pt}

	\begin{minipage}{\textwidth}
	\begin{center}
			\bf{Bouncing on the top side of the square}
	\end{center}
	\makebox[\linewidth]{	
		\includegraphics[width=3.9in]{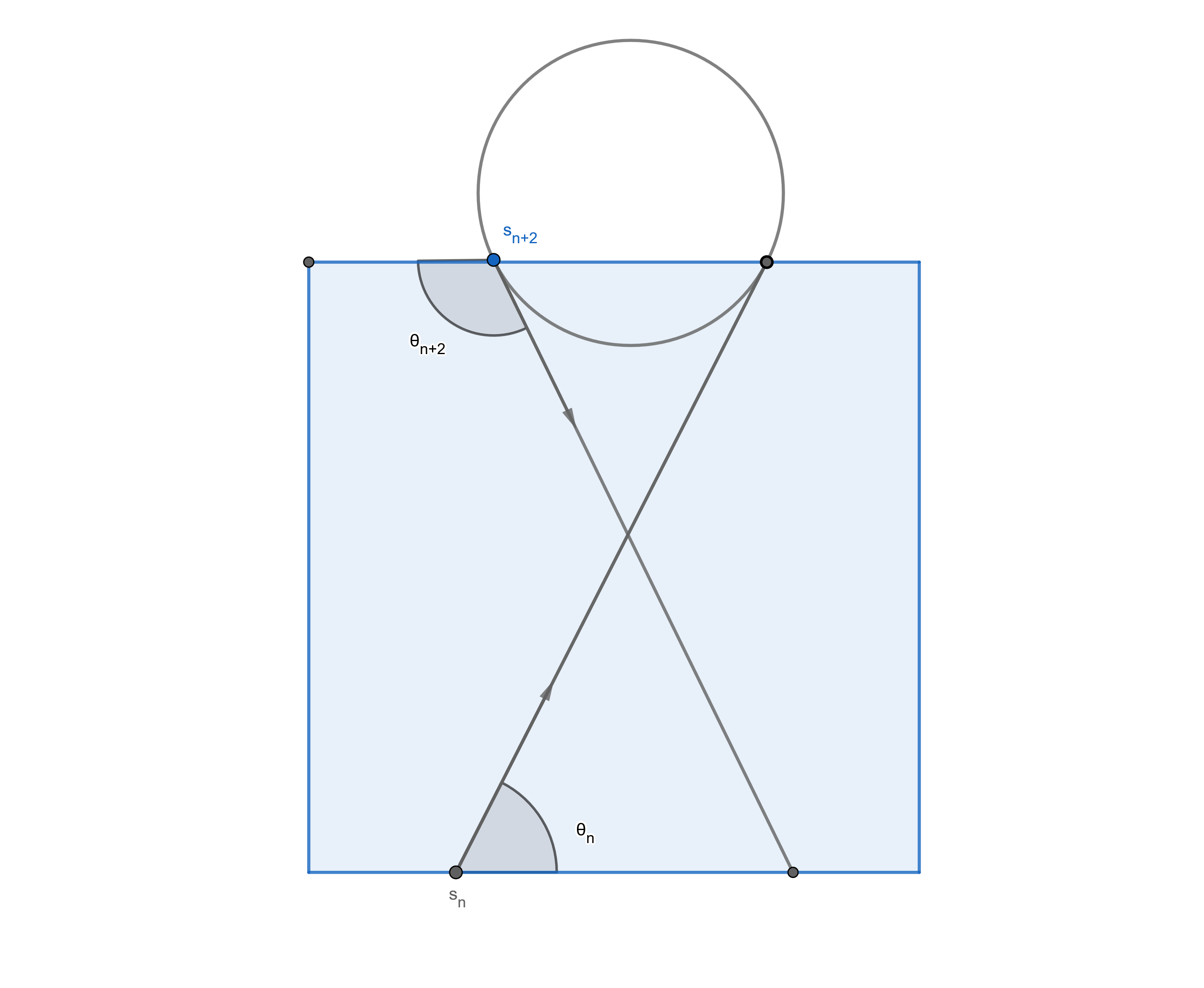}}
		\vspace{-20pt}
		\captionof{figure}{Bouncing on the top side}\label{topside}
	\begin{equation*}
		s_{n+1} = 1-s_n - \cot\theta, \quad
		s_{n+2} = s_{n+1} + \frac{2}{B}\sin\theta,\quad
		\theta_{n+2} = \pi-\theta_n.
	\end{equation*}
	\end{minipage}

	\begin{minipage}{\linewidth}
		\begin{center}
			\bf{Bouncing on the left side of the square}
		\end{center}
	\makebox[\linewidth]{	
		\includegraphics[width=4in]{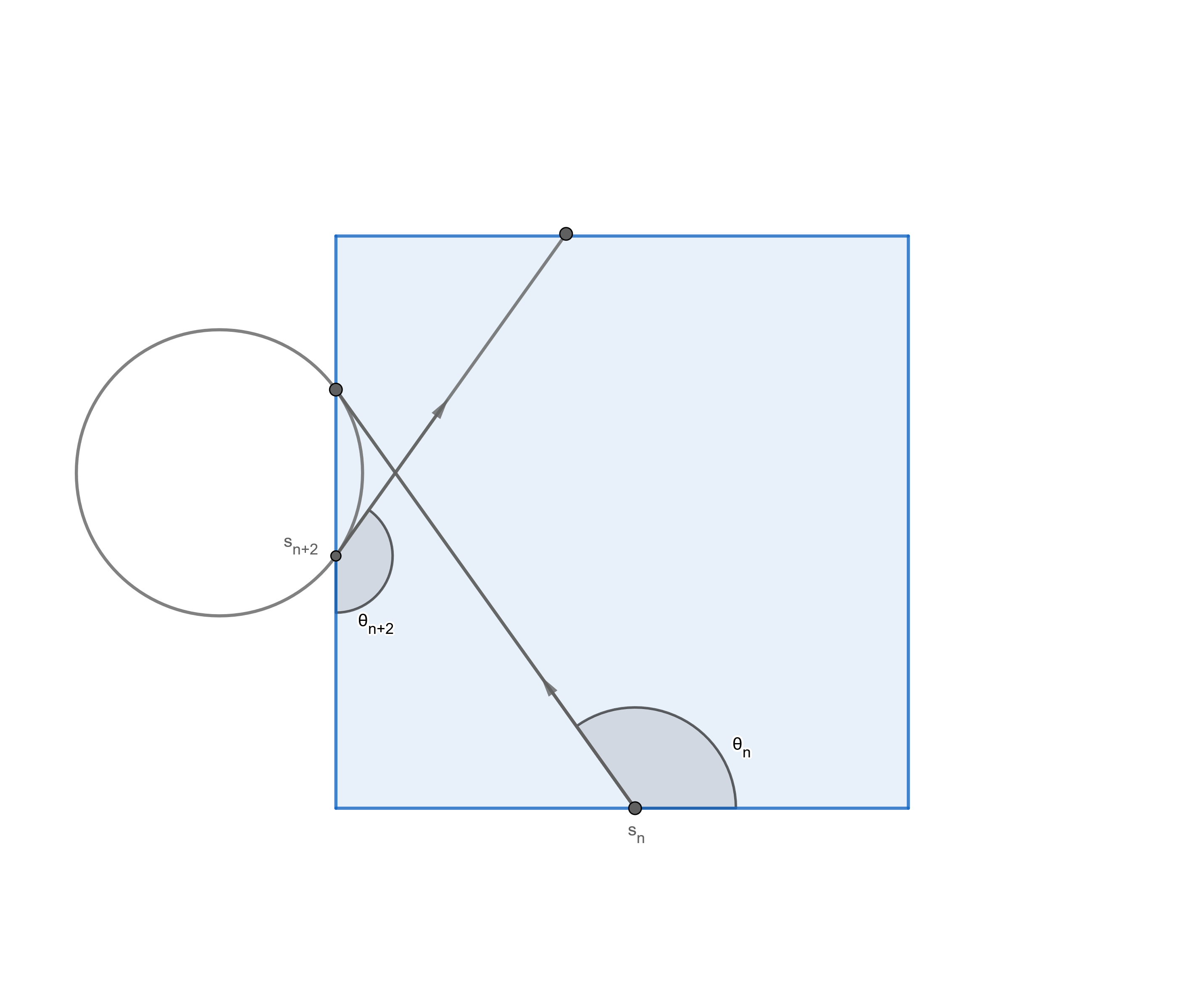}}
		\vspace{-20pt}
		\captionof{figure}{Bouncing on the left side}\label{leftside}
	\begin{equation*}
		s_{n+1} = 1-s_n\tan\theta, \quad
		s_{n+2} = s_{n+1} + \frac{2}{B}\cos\theta, \quad
		\theta_{n+2} = \frac{3\pi}{2}-\theta_n.
	\end{equation*}

	\end{minipage}

	\subsection{Rational slope ($\boldsymbol{\tan\theta}\in\boldsymbol{\mathbb{Q}}$)}
		
	When the slope is a rational number $\frac{p}{q}$ with $p$ and $q$ relative primes, for classical billiards the orbit is periodic. For our inverse magnetic billiard we'd like the nearby periodic orbits, to understand the problem we want conditions where the particle never turns a corner.\\\ \\
	
	Not turning a corner is the same as saying that the unfolded trajectory will cross the same sides that the classical one, a straight line in $\mathbb{Z}\times \mathbb{Z}$, the will never be at opposite sides of a vertex.\\
	
	Since the slope is $\frac{p}{q}$ we will have $2(p+q)$ bounces, for us these bounces are actually exits and returns from and to the unit square. The figure \ref{23loop} shows the case of slope $2/3$, where we have $2(2+3) =10$ bounces in the classical billiard.

	\begin{minipage}{\linewidth}
	\makebox[\linewidth]{	
		\includegraphics[width=4in]{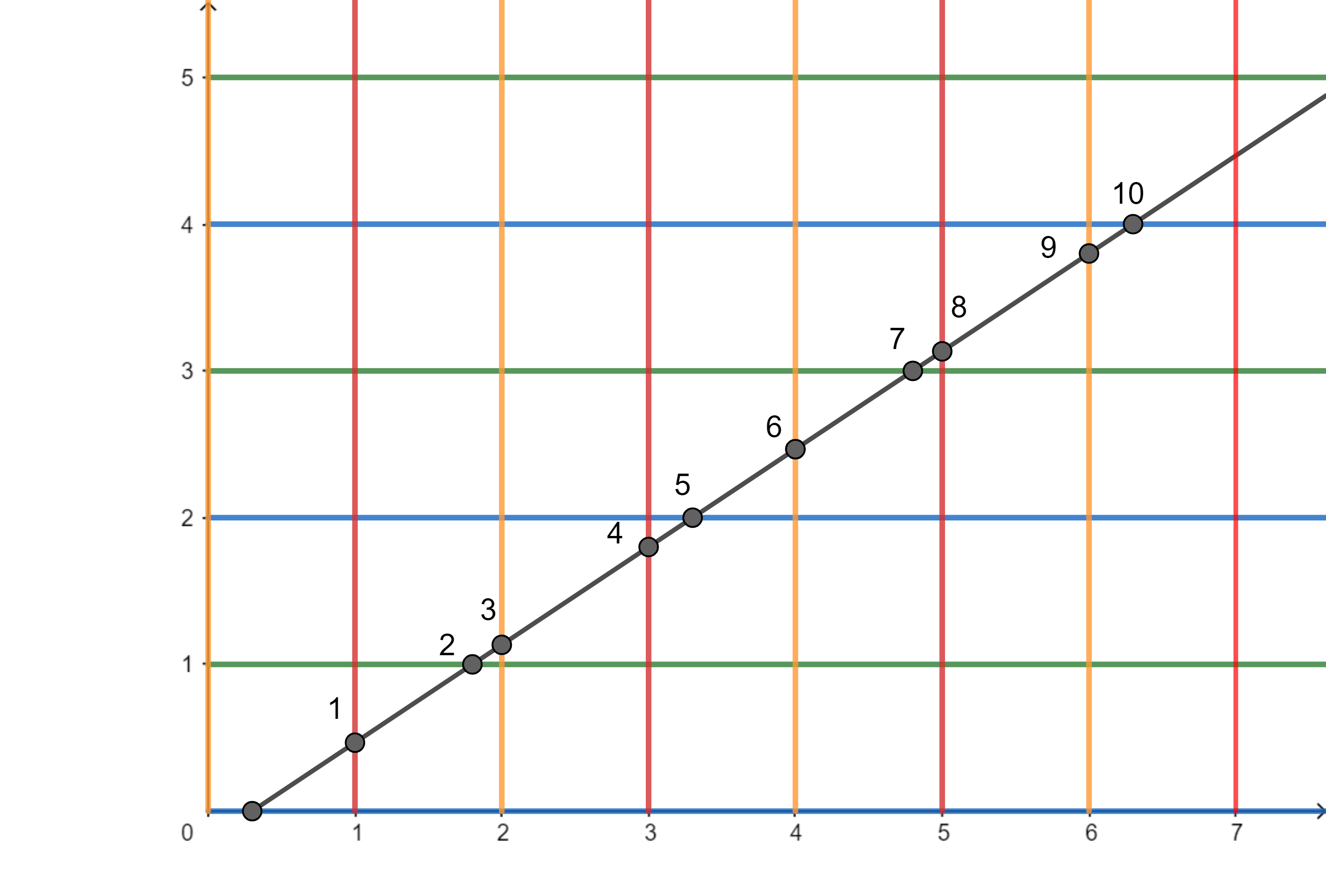}}
		\captionof{figure}{Case of slope $\frac{p}{q} = \frac{2}{3}$}\label{23loop}
	\end{minipage}

	\vspace{10pt}	
	
	The unfolded lattice for the magnetic billiard with these conditions is as shown in the figure \ref{mag23loop}
	
	\begin{minipage}{\linewidth}
		\makebox[\linewidth]{	
			\includegraphics[width=4in]{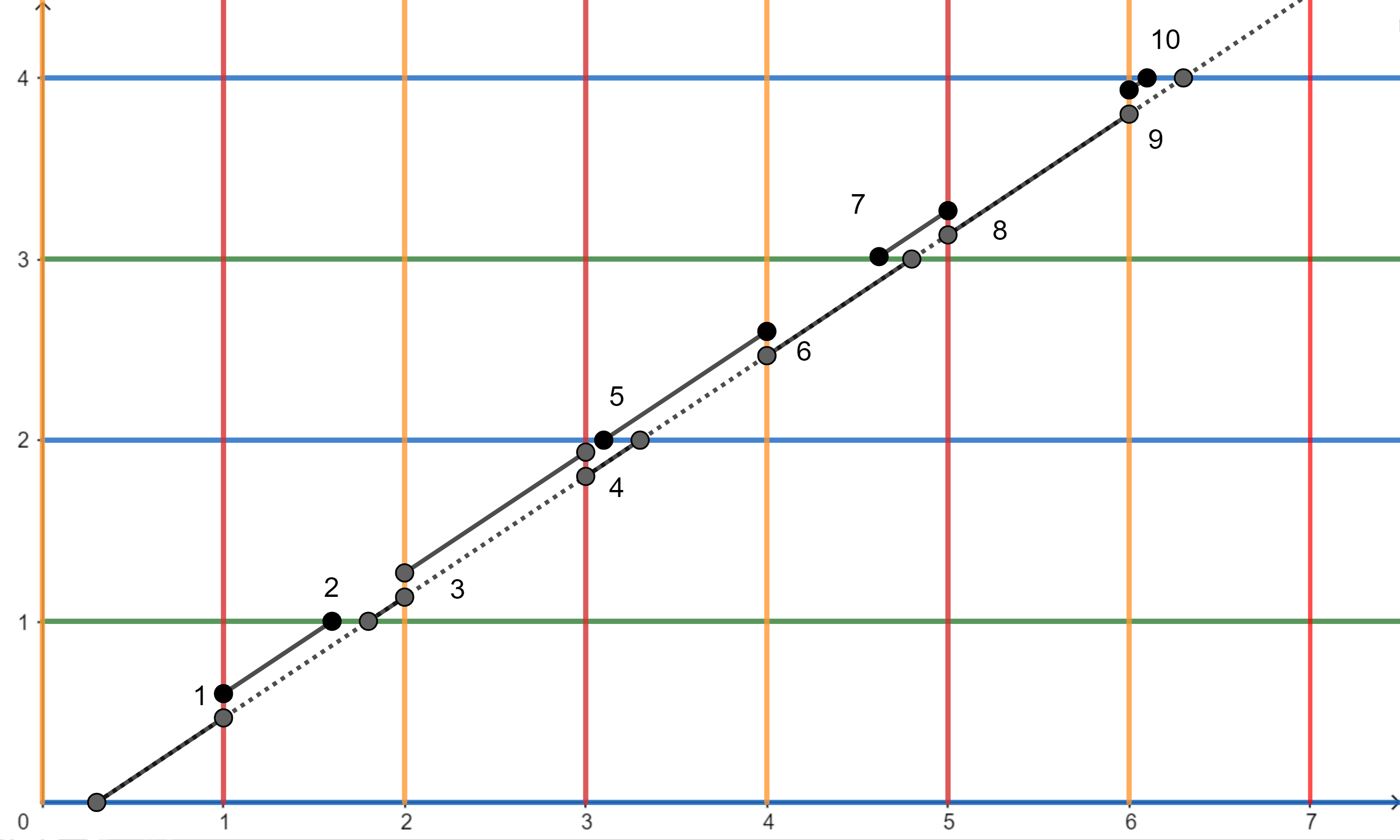}}
		\captionof{figure}{Case of slope $\frac{p}{q} = \frac{2}{3}$}\label{mag23loop}
	\end{minipage}
	
	\vspace{10pt}
	
	\textit{Lemma: With initial conditions $(s_0,\theta_0)$ with $\tan \theta_0 = \frac{p}{q} \in \mathbb{Q}$, $s_0 \neq \mathbb{Z} p$ and $2/B < \min \{ |s_0 - k/q | , |s_0p/q - k/p| : k \in \mathbb{Z}\}$, then the orbit is periodic.} 
	\\\ \\
	
	The condition $s_0 \neq \mathbb{Z} p$ deletes the cases where the orbit hits corners. We can take $B> 2 / \min \{ |s_0\frac{p}{q} - \frac{1}{q}|, |s_0 - \frac{1}{p}|   \}$ so we actually insure that the radius of the circle is small enough to stay in the same side of the square. Since $\frac{2}{B}$ is the diameter of the circle, with this condition we are securing that the exiting point $(s_1,\theta_1)$ is at least at one diameter distance from all the corners. This choice of $B$ makes the diameter of the circle smaller than all the distances between the orbit and the points in the lattice $\mathbb{Z}\times \mathbb{Z}$. \\
		
	The closest that a straight line that starts at the origin with slope $\frac{p}{q}$ gets to a point in the lattice is $\min \{ |k/q| , |k/p| : k \in \mathbb{Z}\}$. Then a line that starts at $(s_0,0)$ will get shifted at the intersection points (with the lines in the lattice) by $s_0$ to the right on the horizontal lines, and by $s_0\frac{p}{q}$ downwards in the vertical lines. This means that we can get as close as we want to the points in the lattice depending on $s_0$. Once you fix $s_0$, the closest you get is $\min \{ |s_0 - k/q | , |s_0p/q - k/p| : k \in \mathbb{Z}\} $.\\\
		
	In our setting we will have $2(p+q)$ shiftings of the classical trajectory, these shifting  can be $\epsilon_1=\frac{2}{B}\cos\theta_n$ or $\epsilon_2=\frac{2}{B}\sin\theta_n$ depending of the side of bounce: $\epsilon_1$ for a vertical side (red or orange in the grid) and $\epsilon_2$ for horizontal sides (blue and green on the grid). These numbers correspond to the formulas in the previous section.\\

	In the unfolded lattice $\mathbb{Z}\times \mathbb{Z}$ the trajectory will cross $2q$ vertical lines and $2p$ horizontal lines. $p$ of those crossing correspond to the lower side of the square, $p$ to upper side, $q$ to the left side and $q$ to the right side. In our notation this means $2p$ shifts of $\epsilon_2$ units horizontally and $2q$ shifts of $\epsilon_1$ units vertically.\\
	
	Since the rotation is counter clockwise, the shifting on the lower side is always to the right, on the right side is upwards, on the top side is to the left and on the left side is downwards. Since the shifting is towards opposite directions on parallel sides they will cancel each other. With this we have that $\displaystyle (s_0,\theta_0) = (s_{2(p+q)},\theta{2(p+q)})$, this proves the lemma.\\

	The figure shows the first bounces for the example $p/q=2/3$. Every shift has another one that cancels it out, after $2(p+q)$ bounces we are at the initial conditions again.\\
	
	\begin{minipage}{\linewidth}
		\makebox[\linewidth]{	
			\includegraphics[width=4in]{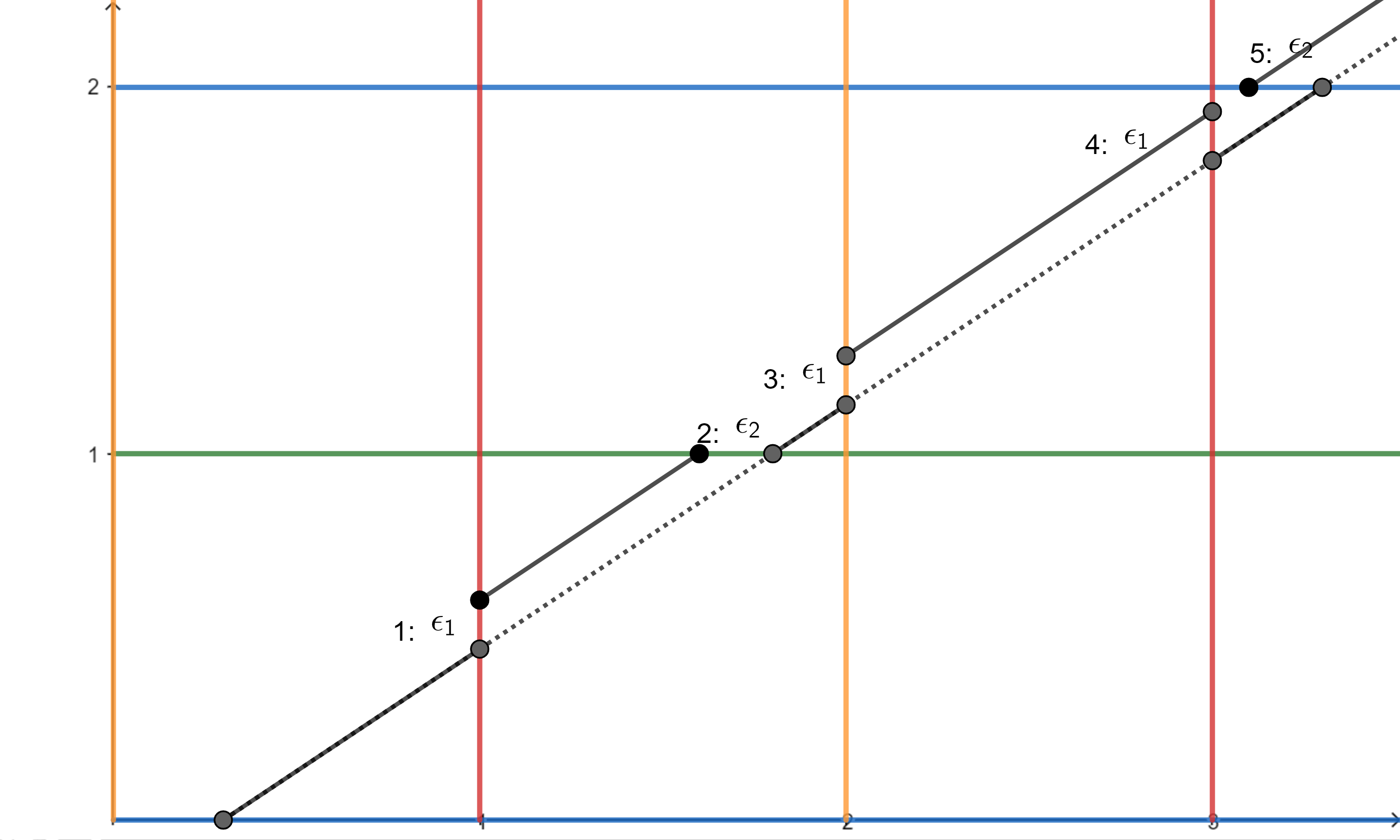}}
		\captionof{figure}{Case of rational slope $\frac{p}{q}$}
	\end{minipage}
	
	\vspace{10pt}

	\subsubsection{Irrational slope}
	
	This case is impossible to have with all bounces on the same exiting side. The particle must turn a corner: the same argument as before works here, the shifting of a classical dense orbit in our tessellation will give us a turn of a corner.\\
	
	Examples of this case will be shown in the next section.
	
	\subsection{Lemma,	$B>1$.}
	
	\textit{Lemma}:	For any radius smaller than 1 (i.e. for any magnetic field with magnitude  bigger than 1), there exist a periodic orbit.\\\ \\
	\textit{Proof}: Pick the initial condition $(r,\pi/2)$.\\

\section{Numerics: Turning corners, turning chaotic}

	In this section we give numerical evidence that our billiard map tends to be chaotic.\\
	After a turn of a corner the set of $\theta$ values is dense, the slope can be rational or irrational. If rational, we are not in the conditions of our lemma, we have violated the condition of being away from a corner by the minimum established distance.\\
	
	\subsection{$B>1$}
	
	In this case we are demanding the radius of the circle to be small enough to go around at most one corner, so the entering side is the same or adjacent to the exiting side.\\
	
	We calculate the map $F:\Sigma^2 \to \Sigma^2$ with different initial conditions (on the Bottom side) and radius. The pictures show $F^n(s_0,\theta_0) = (s_{2n},\theta_{2n})$ for different numbers of iterations. Here we use the notation $u=\cos\theta$.\\
	
	\subsubsection{Radius $0.02$.} Let's look to a periodic orbit first:
	
	\begin{minipage}{\linewidth}
	\makebox[\linewidth]{	
		\includegraphics[width=4in]{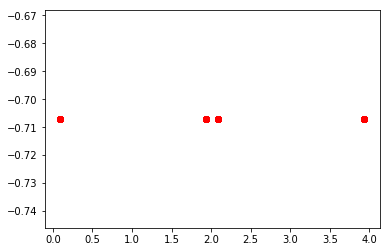}}
		\captionof{figure}{$r=0.02, s = 0.9, \theta = \pi/4$}
	\end{minipage}

	\vspace{10pt}

	This one is similar to figure \ref{slope1nocorners} with slope $1$.\\
	Now four different orbits (4 different initial conditions):
	
	\begin{minipage}{\linewidth}
	\makebox[\linewidth]{	
		\includegraphics[width=4in]{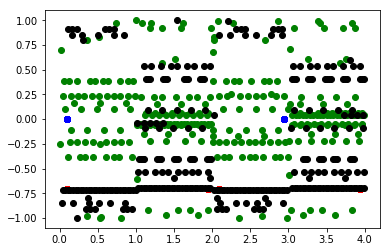}}
		\captionof{figure}{$r=002; (s,u)=(0.1,0),(0.3,0.8), (0.5,-0.5), (0.09,\sqrt{2}/2)$}
	\end{minipage}

	\vspace{10pt}

	Now with more iterations ($1000$) was done with several initial conditions with the same result:
		
	\begin{minipage}{\linewidth}
	\makebox[\linewidth]{	
		\includegraphics[width=4in]{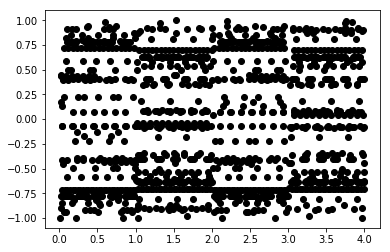}}
		\captionof{figure}{$r=0.02, s=0.001, u=-0.99$ REP = 1300}
	\end{minipage}

	\vspace{10pt}

	Going for $10000$ bounces
	
	\begin{minipage}{\linewidth}
		\makebox[\linewidth]{	
			\includegraphics[width=4in]{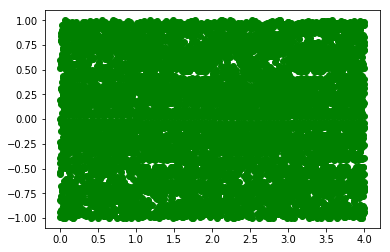}}
		\captionof{figure}{$r=0.02, s=0.001, u=-0.99$ REP = 10000}
	\end{minipage}
		
	\subsubsection{Radius $0.49$}
	
	Now we try with a different radius. Here we encounter a piece-wise linear behavior. If you increase the radius, the pattern becomes clearer.
	
	\begin{minipage}{\linewidth}
	\makebox[\linewidth]{	
		\includegraphics[width=4in]{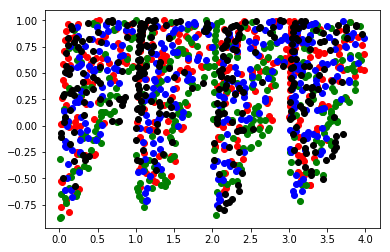}}
		\captionof{figure}{r=0.49, s=0.1 u=0, s=0.2 u=0.8, s=0.5 u=-0.5, s=0.09 u=0.45}
	\end{minipage}

	\vspace{10pt}
	
	\begin{minipage}{\linewidth}
	\makebox[\linewidth]{	
		\includegraphics[width=4in]{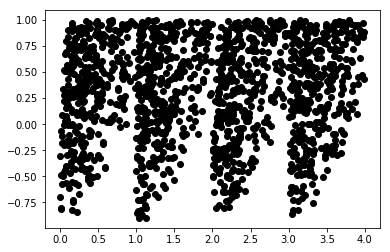}}
		\captionof{figure}{r=0.49, s=0.3, u=0.8}
	\end{minipage}
	
	\vspace{10pt}
	
	\subsubsection{Radius $0.85$}
	
	More clear that there is a piecewise linear behavior.
	
	\begin{minipage}{\linewidth}
	\makebox[\linewidth]{	
		\includegraphics[width=4in]{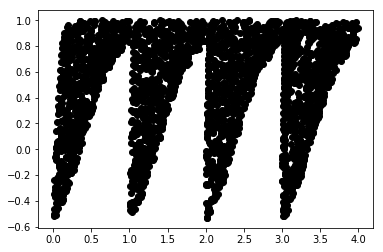}}
		\captionof{figure}{r=0.85, s=0.33, u=0.8, 3000 repetitions}
	\end{minipage}

\section{Conclusion}

Inverse magnetic billiard dynamics limits to regular billiards on the square, but turns out to be much more intricate and complicated than regular billiards. Numerical experiments suggest that the dynamics is ergodic for strong magnetic fields, while it has  unexplained sawtooth type patterns in the intermediate range of fields.

As far as our original questions, listed at the the end of section 1, we have answered  the first three, those regarding periodic orbits, affirmatively.  We are unable to classify periodic orbits that turn corners, although by symmetry these exist. We also found dense orbits in both spaces $\Omega$ and $\Sigma^2$, we didn't classify these. The remaining questions remain open.

\section{}

\vspace{20pt}


\begin{thebibliography}{20}
	
	\bibitem{Arnold} V. Arnold, \textit{Mathematical methods of classical mechanics}, Springer, 1989.
		
	\bibitem{Birkhoff} G.D. Birkhoff, \textit{Dynamical Systems}, American Mathematical Society / Providence, RI, American Mathematical Society, 1927.

	
	\bibitem{BerglundHunz} N. Berglund and H. Kunz, \textit{Integrability and ergodicity of classical billiards in a magnetic field}, Journal of Statistical Physics 83, no. 1-2, 81–126, 1996.
	
	\bibitem{ChernovMarkarian} N. Chernov and Ro. Markarian, \textit{Chaotic Billiards}, American Mathematical Soc., 2006.
	
	\bibitem{DatserisHupeFleischmann}	G. Datseris, L. Hupe,and R. Fleischmann, \textit{Estimating Lyapunov exponents in billiards}, Chaos: An Interdisciplinary Journal of Nonlinear Science, 29 no.9, p.093115, 2019.
	
	\bibitem{Gasiorek} S. Gasiorek, {\it On the Dynamics of Inverse Magnetic Billiards}, Ph.D. thesis, University of
	California Santa Cruz, 2019
		
	\bibitem{RobnikBerry} M. Robnik and M. V. Berry, \textit{Classical billiards in magnetic fields}, J. Phys. A: Math. Gen. 18, no. 9, 1361–1378, 1985.
	
	\bibitem{Robnik} M. Robnik, \textit{Regular and chaotic billiard dynamics in magnetic fields}, Nonlinear Phenomena and Chaos 1, 303–330, 1986.
	
	\bibitem{Tabachnikov95} S. Tabachnikov, \textit{Billiards}, Soc. Math. France, 1995.
	
	\bibitem{Tabachnikov93} S. Tabachnikov. \textit{Outer billiards}, Russ. Math. Surv., 48, 75-102, 1993.
	
	\bibitem{Tabachnikov04} S. Tabachnikov, \textit{Remarks on magnetic flows and magnetic billiards, Finsler metrics and a magnetic analog of Hilbert’s fourth problem}, Modern dynamical systems and applications, pp.233-250, 2004.
		
	\bibitem{Tasnadi97} T. Tasnádi, \textit{Hard Chaos in Magnetic Billiards (On the Euclidean Plane)}, Communications in Mathematical Physics 187, no. 3, 597–621, 1997.

	\bibitem{Tasnadi96} T. Tasnadi, \textit{The behavior of nearby trajectories in magnetic billiards}, J. Math. Phys. 37, 5577-5598, 1996.

	\bibitem{VorosTasnadiCsertiPollner} Z. Vörös, T. Tasnádi, J. Cserti, and P. Pollner, \textit{Tunable Lyapunov exponent in inverse magnetic billiards}, Physical Review E 67,	no. 6, 065202, 2003.
	
	
\end{thebibliography}
\end{document}